\documentclass[11pt]{amsart}
\usepackage{amsfonts,amsmath,amssymb,url,verbatim,multirow}
\usepackage{algorithm,algorithmic}
\newcommand\E{\ensuremath{\mathbb{E}}}
\newcommand\bT{\ensuremath{\mathbb{T}}}

\newcommand\R{\ensuremath{\mathbb{R}}}

\newcommand\PP{\ensuremath{\mathbb{P}}}

\newcommand\F{\ensuremath{\mathcal{F}}}
\newcommand\U{\ensuremath{\mathcal{U}}}
\newcommand\G{\ensuremath{\mathcal{G}}}

\newcommand\Z{\ensuremath{\mathcal{Z}}}

\newcommand\vtau{\ensuremath{\vec{\tau}}}

\newcommand\vp{\ensuremath{\vec{\pi}}}

\newcommand\vze{\ensuremath{\vec{\zeta}}}
\newcommand\valpha{\ensuremath{\vec{\alpha}}}
\newcommand\utau{\ensuremath{\underline{\tau}}}

\newcommand\hq{\ensuremath{\hat{q}}}
\newcommand\vth{\ensuremath{\vartheta}}
\newcommand\tJ{\ensuremath{\tilde{J}}}
\newcommand\hJ{\ensuremath{V}}

\newcommand\sD{\ensuremath{\mathfrak{D}}}

\newcommand\vxi{\vec{u}}

\newcommand\dt{\;} 
\newcommand\eps{\ensuremath{\epsilon}}

\DeclareMathOperator*{\argmin}{arg\, min}

\newtheorem {thm}{Theorem}[section]
\newtheorem {prop}{Proposition}[section]
\newtheorem {lemm}[thm]{Lemma}
\newtheorem {deff}{Definition}[section]
\newtheorem {cor}[thm]{Corollary}

\theoremstyle{remark}
\newtheorem {rem}{Remark}[section]

\usepackage[margin=1in]{geometry}

\usepackage[pdftex]{graphicx}

\title[Stochastic Switching Games]{Stochastic Switching Games and Duopolistic Competition in Emissions Markets}

\author{Michael Ludkovski} \address[M. Ludkovski]{Department of Statistics and Applied Probability,
University of California, Santa Barbara, CA 93106-3110, USA} \email{ludkovski@pstat.ucsb.edu,
www.pstat.ucsb.edu/faculty/ludkovski}

\begin{document}
\begin{abstract}
We study optimal behavior of energy producers under a $CO_2$ emission abatement program. We focus on a two-player discrete-time model where each producer is sequentially optimizing her emission and production schedules. The game-theoretic aspect is captured through a reduced-form price-impact model for the $C O_2$ allowance price. Such duopolistic competition results in a new type of a non-zero-sum stochastic switching game on finite horizon. Existence of game Nash equilibria is established through generalization to randomized switching strategies. No uniqueness is possible and we therefore consider a variety of correlated equilibrium mechanisms. We prove existence of correlated equilibrium points in switching games and give a recursive description of equilibrium game values. A simulation-based algorithm to solve for the game values is constructed and a numerical example is presented.

%
\end{abstract}


\maketitle

\section{Introduction}
In this paper we study a new class of non-zero-sum stochastic \emph{switching games} with continuous state-space. Such games have natural applications in economics and finance, in particular for describing oligopolistic competition between large commodity producers. Our analysis is motivated by the $C O_2$ cap-and-trade markets and provides new quantitative insight into the game-theoretic aspects of these schemes.

Switching game are a special class of dynamic non-zero-sum state-space games. They are characterized by a finite number of system states $\vxi$, jointly selected by the players. The players dynamically react to actions of other players and the evolution of state variables, represented as controlled stochastic processes, by strategically modifying the system state. Our contribution is a first rigorous probabilistic analysis of switching games. Because multiple game Nash equilibria are possible in our model, we propose to apply the wider concept of correlated equilibria. Correlated equilibria give a clear financial mechanism for stepwise equilibrium selection. Our key result is the construction of correlated equilibria in switching games in Section \ref{sec:multiple-stopping-game}. The resulting representation in Theorem \ref{thm:iter-game} of switching games in terms of a recursive sequence of stopping games leads to a constructive characterization of equilibrium strategies. Namely, we prove the analogue of the dynamic programming equation for the game values which enables numerical solution through backward recursion. Thus, the complexity of switching games is only slightly higher than of regular optimal switching problems.

In terms of existing literature this paper extends two separate strands of research. Work on stochastic zero-sum  \emph{stopping games} dates back to Dynkin \cite{DynkinGame69}. Such Dynkin games were progressively generalized in \cite{Alvarez08,KaratzasCvitanic96,EkstromPeskir08,Ferenstein07,TouziVieille02}. Later extensions also treated special cases of non-zero-sum stopping games, especially the so-called monotone type \cite{HamadeneZhang09,Ohtsubo87,Ohtsubo91}. The key tool of correlated equilibria in stochastic dynamic games was studied by \cite{NowakRaghavan92,RamseySzajowski08,RosenbergSolanVieille,ShmayaSolan,SolanVieille02,SolanVohra01}. We augment these results by explicitly characterizing correlated equilibria in repeated stopping games using the methods of Ramsey and Szajowski \cite{RamseySzajowski08}.
Contemporaneously, the theory of \emph{optimal switching} for a single agent was developed and extensively studied in the past decade, see \cite{CarmonaLudkovskiSwitching,DayanikEgamiSwitching,HamadeneJeanblanc,Pham06}. In Section \ref{sec:multiple-stopping-game} we extend these results to a game setting by showing that at switching game equilibrium each player faces a optimal switching problem with \emph{randomized} controls.

Another contribution of this work is a construction of numerical schemes to compute game-values and equilibrium strategies of switching games. This is achieved by combining backward recursion together with sequential solution of local $2 \times 2$ games. We suggest two approaches, one based on the Markov chain approximation method and a second algorithm that relies on least squares Monte Carlo. The latter is a novel extension of our previous work in \cite{CarmonaLudkovskiSwitching,LudkovskiFlex} and
borrows ideas from standard optimal stopping theory to implement the analogue of the dynamic programming recursion on a set of Monte Carlo simulations.

A significant portion of the paper is dedicated to the application of our model to emissions trading.  With imminent ramping up of $C O_2$-emissions markets around the world (see e.g.~ the Western Climate Initiative in the US and the EU ETS Phase III in Europe both set to start in 2012), it is crucial to understand energy producer behavior under the new frameworks. By design, the carbon allowances will be scarce and market participants will be competing for finite permit resources. Our analysis is a first pass at oligopolistic competition in $C O_2$ markets using game-theoretic methods. We hope it can serve
as a stepping-stone to more sophisticated modeling that addresses market design and comparative statics of our framework.


The rest of the paper is organized as follows. In Section \ref{sec:competition} we define the precise stochastic model representing oligopolistic competition in $CO_2$ markets.
Section \ref{sec:seq-games} constructs the representation of switching games in terms of repeated stopping games and culminates with Theorem \ref{thm:iter-game} that establishes the dynamic programming equations. Section \ref{sec:numerical} describes our numerical solution algorithms and presents a computational example. Finally, Section \ref{sec:extensions} discusses extensions of the model and points directions for future work.

\section{Competitive Equilibrium among Oligopolistic Producers}\label{sec:competition}
In this section we formally define a model for the competitive dynamic equilibrium between producers with market power.

\subsection{Price Dynamics under Cap-and-Trade Emission Schemes}\label{sec:emission-trading}
The carbon allowance markets are intrinsically linked to other energy markets, notably electricity whose production accounts for the bulk of regulated emissions. Furthermore, due to their size, major electricity producers often have the ability to dramatically move carbon prices based on their emission schedules. Hence, to understand equilibria in $C O_2$-markets it is useful to consider them from the point of view of \emph{large traders}.
The model proposed below captures this phenomenon for a joint electricity-carbon market.

To fix ideas, start with a filtered probability space $(\Omega, \mathcal{H}, (\F_t), \PP)$, $t\in \bT \triangleq \{0,1,\ldots,T\}$. The terminal time $T$ is the expiration date of the current vintage of permits. We consider two heterogenous producers (henceforth termed players) who each produce commodity $P$ (electricity) and consume commodity $X$ (carbon allowances). These two producers generate ``dirty'' electricity from e.g.~coal or gas and can be viewed as representative agents of a park of power plants with identical engineering characteristics. All the other participants in the electricity and $CO_2$ markets are not modeled explicitly; rather we postulate that their collective actions induce \emph{stochastic} fluctuations in the prices of $P$ and $X$. We assume that the two players are large traders in the carbon market, but small players on the electricity market. This reflects the fact that the other ``green'' producers (who use nuclear, hydroelectric, renewable, etc.~sources) create a competitive electricity market while remaining passive in the $C O_2$-permits arena.

The electricity price is given by the exogenous discrete-time stochastic process $(P_t)$, for simplicity taken to be one-dimensional,
\begin{align*}
P_{t+1} = G(P_t, \eps^P_t),
\end{align*}
where the innovations $(\eps^P_t)$ are i.i.d.~standard Gaussian.
Our canonical example is the logarithmic Ornstein-Uhlenbeck process which is a log-normal stationary Gaussian process with
\begin{align}\label{eq:p-dynamics}
P_{t+1} = P_t \cdot \exp \left( \kappa_P(\bar{P} - \log P_t)\dt + \sigma_P \eps^P_t \right),\qquad \eps^P_t \sim \mathcal{N}(0,1),
\end{align}
for some positive constants $\kappa_P, \bar{P}, \sigma_P$.

The objective of the producers is to maximize their expected net profits over the planning horizon $T$. The producers' profit is given by their \emph{clean dark spread} \cite{Eydeland} which is defined as the difference between electricity price and the carbon-adjusted production cost. We assume that input fuel costs are fixed, as is often the case for power generators with long-term supply contracts. The strategy of each player is described by a repeated start-up/shut-down option. Namely, if the market conditions are unfavorable, a player can stop production, eliminate $C O_2$ emissions and avoid losses; she can then restart production when the profit spread improves. As a first approximation we assume that these choices of production regimes are binary and denoted as ``off'' (0) and ``on'' (1). Formally, the production schedule of each producer is described by a stochastic process $u_i$, $u_i(t) \in \{0,1\}$, $t \in \mathbb{T}$.
In a single-player model, such timing optionality is known as a {real option} and has been thoroughly investigated since the seminal work of \cite{BrennanSchwartz,DixitPindyck}. Repeated real options have attracted considerable attention recently, see \cite{CarmonaLudkovskiSwitching,HamadeneJeanblanc,LudkovskiFlex}, and others. 

\begin{rem}
An alternative formulation is to take $u_i$ to be continuous, so that the producers can choose emissions levels smoothly.
This would lead to a non-zero-sum stochastic dynamic game. Such models have been extensively studied both in discrete and continuous time, see
~\cite{HamadeneLepeltier00,Nowak03} and references therein. While presenting other formidable technical challenges, the problem of equilibrium selection is less severe with continuous controls thanks to the convexity of value functions.
In this work we focus on the timing flexibility and therefore maintain the discrete control space.
%
%
%
\end{rem}

Let $X_t$ be the permit price at date $t$. Based on above discussion, the actions of each player influence the dynamics of $X_t$. Namely,  conditional on player actions $u_1(t), u_2(t)$, we model $X_t$ as another mean-reverting process with a policy-dependent mean and log-Gaussian increments,
\begin{align}\label{eq:X-dynamics} \left\{ \begin{aligned}
X_{t+1} &= X_t \cdot \exp \left( \kappa_X( f(u_1(t), u_2(t)) - \log X_t)\dt + \sigma_X \eps^X_t \right) \qquad\text{with}\\ f(u_1, u_2) & = \log(\bar{X} + g_1 u_1 + g_2 u_2). \end{aligned}
\right.
\end{align}
The sequence $(\eps^X_t)$ is again Gaussian, with correlation parameter $\rho$ to $(\eps^P_t)$,
i.e.~$\eps^X_t = \rho \eps^P_t + \sqrt{1-\rho^2} \eps^\perp_t$ with $\eps^\perp_t \sim \mathcal{N}(0,1)$ independent of $\eps^P_t$. Rising electricity prices are likely to increase the overall $C O_2$ emission rates and therefore we expect that $P$ and $X$ are positively correlated, $\rho >0$.

\begin{rem}
To motivate \eqref{eq:X-dynamics}, we recall from \cite{FehrHinz08} that in a carbon market ``$X_t = \overline{x} \bar{\PP}\{  C_T > \bar{c} | \F_t\}$,'' where $C_T$ is the cumulative total $CO_2$ emissions on $[0,T]$, $\bar{x}$ is the penalty for going over the allowance limit, $\bar{c}$ is the total amount of allowances allocated and $\bar{\PP}$ is the equilibrium pricing measure. We postulate that $C_T = \sum_{s=0}^{T-1} \left\{ b_1 u_1(s) + b_2 u_2(s)  + \tilde{u}(s) \right\}$, where $b_i u_i(s)$ are the emissions of producer $i$ in period $s$ and $\tilde{u}(s)$ are the emissions by all other market participants. Assuming independent increments (due to external shocks such as weather effects, etc.) in $\tilde{u}(s)$, the dynamics \eqref{eq:X-dynamics} follow, with some complicated and time-dependent functions $f(t,\cdot)$ and volatility $\sigma_X(t)$. In \eqref{eq:X-dynamics} we give a simplified or reduced-form version of this description to capture the temporal feedback between $u_i$'s and $X_t$.
If the supply curve for the $C O_2$  allowances is convex, then the price impact $f(u_1,u_2)$ in \eqref{eq:X-dynamics} would be nonlinear and further magnify the competitive effects.
\end{rem}

\subsection{Optimization Objective}
We assume that the producers have zero allowance allocations and cannot bank allowances; therefore they must purchase the requisite allowances at each stage of the game. The total P\&L of the players then consists of
(i) revenue from selling electricity, minus the (ii) cost of buying emission allowances, as well as (iii) operational costs due to adopted strategy $u_i$. An important case of operational costs are \emph{fixed switching costs} $K_{\{i,j_1,j_2\}}$ that are paid each time the production regime of agent $i$ is changed from $j_1$ to $j_2$ and corresponding to the ramping-up/winding-down costs associated with the electricity turbines \cite{CarmonaLudkovskiSwitching,Eydeland}. We postulate $K_{\{i,j_1,j_1\}} = 0 \;\forall j_1$ and the triangle inequality $K_{\{i,j,\ell\}} \le K_{\{i,j,k\}} + K_{\{i,k,\ell\}}$ for all $j,k,\ell$.

Let $\PP^{\vxi}$ be the law of $(P_t, X_t)$ given a strategy pair $\vxi \equiv (u_1(t), u_2(t))_{t=0}^T$.
The expected cumulative net profit of producer $i$ starting with $P_s=p$, $X_s = x$ and initial production regime $\vze \in \{0,1\}^2$ is
\begin{align}\label{eq:expected-profit}
V_i(s, p, x, \vze; \vxi) & \triangleq \E^{\vxi} \left[ \sum_{t=s}^{T-1} \left\{(a_i P_t - b_i X_t - c_i) u_i(t) \, \dt - K_{\{i, u_i(t-), u_i(t)\}}  \right\} \Big| P_s = p, X_s = x, \vxi(s-) = \vze\right].
\end{align}
The constants
 $a_i, b_i, c_i, K_{\{i,j_1,j_2\}}$, $i=1,2$, represent the maximum quantity of electricity produced by the facility in one period, the amount of corresponding $C O_2$ allowances needed,
 fixed production costs and switching costs, respectively. Due to switching costs, current production regime is also a state variable.
 Below, the theorems on existence of equilibria in stochastic games require bounded payoffs; therefore we assume that profits are truncated from above at some large positive level.


%




In the duopoly setting, while each player aims to maximize her own profits, the competitor actions will also affect her decisions.
Indeed, emissions today shrink remaining permit supplies and tend to increase \emph{future} $C O_2$ prices. Therefore, if player 1 is emitting, player's 2 expected future profits are reduced.
%
Overall, the producers are facing a stochastic game where actions correspond to the latest choice of production regime by each player and payoffs are a function of the exogenous $P_t$ and the partly controlled $X_t$. Since the game is stochastic and multi-period with Markov state variables, we restrict our attention to Markovian (feedback) equilibria. Our main task for the remainder of the paper is to characterize such game equilibria and then compute the corresponding game value functions (i.e.~expected profits) $V_i$ and equilibrium emission schedules $(u_1^*, u_2^*)$.

\subsection{Randomized Emission Schedules}\label{sec:randomized-strat}
The strategies $u_i$ may be mixed or randomized, i.e.\ $u_i(t)$ is not necessarily adapted to the market filtration $\F$. However, we also assume a full-information setting, whereby the emission schedules of each agent are publicly known after the fact. Accordingly, the market observables are $\F_t = \sigma(X_0, P_0, u_1(0), u_2(0), \ldots, u_1(t-1), u_2(t-1), X_t, P_t)$, the filtration generated by the price histories and past actions. An $\F$-randomized emission strategy is a pair $(u_i(t), \G^i(t))$ where $\G^i$ is an independent enlargement of the filtration $\F$ (i.e.~$\PP( A | \F_t) = \PP( A | \tilde{\G^i_t})$ for all $A \in \F_t$) and $u_i$ is $\G^i$-adapted.
Let
\begin{align}
p_i(t) &\triangleq  \PP( u_i(t) = 1 | \F_t)
\end{align}
denote the stage-$t$ probability that the control will be `on', given observable information so far. If $p_i(t) \in \{0,1\}$ then the strategy is \emph{pure} at stage $t$, otherwise it is \emph{mixed} and can be represented via a randomization parameter $\eta_i(t)$ as
\begin{align}\label{eq:eta-sequence}
u_i(t) = 1_{\{ \eta_i(t) \le p_i(t)\}}, \qquad \eta_i(t) \sim Unif(0,1), \eta_i(t) \perp \F_t.
\end{align}
The full mixed strategy is the vector $\vp_i(t) \equiv( 1-p_i(t), p_i(t))$ belonging to  the 2-simplex $\vp_i(t) \in \Delta_2 \triangleq \{ (\pi^0, \pi^1): \pi^j \ge 0, \pi^0 + \pi^1 = 1\}$. The joint action is given by the \emph{strategy profile} $\vp(t)$, with $\pi_i^j(t)$ denoting the probability that player $i$ emits at level $j$. 


The set $\U_i$ of admissible production schedules for player $i$ consists of  $\F$-randomized $\{0,1\}$-valued processes and can be canonically identified with an $\F$-adapted process $(p_i(t))$, $0 \le p_i(t) \le 1$ and independent sequence $\eta_i(t)$ as in \eqref{eq:eta-sequence}. Let $\sD^i(t)$ denote the set of $\G^i$-stopping times bigger than $t$. Because $u_i(t) \in \{0,1\}$, we have a one-to-one correspondence between admissible $u_i$'s and sequences $(\tau^u_k)_{k=1}^\infty$ satisfying $\tau^u_{k+1} \in \sD^i(\tau^{u}_k)$, 
\begin{align}\label{eq:xi-vs-tau}
u_i(t) =  \sum_{k=0}^T u_i(0) 1_{[\tau^u_{2k}, \tau^u_{2k+1})} + (1-u_i(0)) 1_{[\tau^u_{2k+1}, \tau^u_{2k+2})}, \qquad \tau^u_0 = 0.
\end{align}
The switching times $\tau^u_k$ encode the times of production regime shifts defined by $u_i$. The representation \eqref{eq:xi-vs-tau} holds because at most one regime switch can be made by each player at any given stage. Indeed, multiple simultaneous regime switches by the same producer are strongly sub-optimal if $K_{\{i,j_1,j_2\}}>0$ and weakly suboptimal otherwise. A $\G^i$-adapted stopping time $\tau$ can also be viewed as a randomized $\F$-stopping time, via its conditional stopping probabilities $p_t = \PP( \tau = t | \tau > t-1, \F_t)$, namely
\begin{align}\label{eq:randomized-tau}
\tau(p) \triangleq \inf \{ t : \eta_t \le p_t \}, \qquad \eta_t \sim Unif(0,1) \perp \F_t.
\end{align}
When $p_t \in \{0, 1\}$ for all $t$, we are back in the case of regular $\F$-stopping times.

If $\G^1_t \cap \G^2_t = \F_t$ then the randomizations of the two players are independent. Alternatively, correlated decision making can be introduced by making the $\eta_i$'s in \eqref{eq:eta-sequence} dependent. Let $\gamma$ be an $\F$-adapted stochastic process taking values in
$\Delta_4$. Following \cite{RamseySzajowski08} we interpret $\gamma(t)$ as a \emph{weak (stepwise) communication device}, with $\gamma_{ij}(t)$ specifying the probability that player 1 takes action $i \in \{0,1\}$ and player 2 applies action $j \in \{0, 1\}$,
$$
\gamma_{ij}(t) \triangleq \PP( u_1(t) = i, u_2(t) = j | \F_t).
$$
The correlation is implemented via a third party that directs the players to implement a particular action pair through private signals. Namely, the players receive signals
\begin{align}\label{def:signals}\mu_1(t, \gamma) = 1_{\{ \gamma_{10}(t) + \gamma_{11}(t) < \eta(t) \}} \qquad \text{and} \qquad \mu_2(t, \gamma) = 1_{\{ \gamma_{01}(t) + \gamma_{11}(t) < \eta(t)\}},
\end{align}
where $\eta \sim Unif(0,1)$ is only observed by the third party. Setting $u_i(t) = \mu_i(t, \gamma)$, the resulting
strategy profile is denoted $\vxi(t,\gamma)=(\vp_1(t, \gamma), \vp_2(t, \gamma))$ and has dependent marginals and joint distribution $\gamma$. Conditional on the signal at date $t$, an agent can impute the strategy of the other player by e.g.\ $\vp_2(t, \gamma)\big|_{\mu_1(t, \gamma)=0} = \bigl(\frac{\gamma_{0 0}(t)}{\gamma_{00}(t)+\gamma_{01}(t)}, \frac{\gamma_{01}(t)}{\gamma_{00}(t)+\gamma_{01}(t)} \bigr)$. With $\gamma$ in place, the space of randomized strategies is now adjusted to $\U_i(\gamma) = \{ (u_i, \G^i) \in \U_i \;\text{ such that } \G^i_t \supseteq \F_t \vee \sigma(\mu_i(t))$.

\subsection{Correlated Equilibria in $CO_2$ markets}\label{sec:eqm-select}
The introduced correlation mechanism $\gamma$ can be used to define correlated equilibrium points (CEP) in the $C O_2$ emissions duopoly game. To motivate the need for such mechanisms, we observe that intuitively the dynamic switching game is a sequence of one-period bimatrix games. At each stage $t$, the control $u_i(t) \in \{0,1\}$ of each player $i \in \{1,2\}$ is simply ``on/off'', leading to the classic $2 \times 2$ game. From a dynamic point of view, the relevant payoff to the players at stage $t$ is then the sum of the current clean spread and the {continuation value} that corresponds to the game value that can be realized in the future by the respective player contingent on current state of the world. In our repeated game setting, the players must \emph{a priori} agree on how to implement future equilibria, otherwise the computation of continuation values would not be possible. Hence, to have a well-defined switching game value, we need existence-uniqueness of equilibria in the one-period sub-games.

Accordingly, we briefly recall the structure of 2-by-2 one-shot game. Consider the $2\times 2$ game  $H$ with normal form \begin{align}\label{eq:H-def}
H = \begin{pmatrix} (z_1^{00}, z_2^{00}) & (z_1^{01}, z_2^{01}) \\ (z_1^{10}, z_2^{10}) & (z_1^{11}, z_2^{11}) \end{pmatrix},
\end{align}
where the rows of $H$ are chosen by player 1, and the columns by player 2. 
A strategy profile $(\vp^*_1, \vp^*_2)$ is a \emph{Nash equilibrium point} (NEP) of $H$ if we have
$$ \sum_{j,k} \pi^{*,j}_1 \pi^{*,k}_2 z_1^{jk}  = \sup_{\vp_1 \in \Delta_2} \sum_{j,k} \pi^j_1 \pi^{*,k}_2 z_1^{jk}, \quad\text{and}\quad
\sum_{j,k} \pi^{*,j}_1 \pi^{*,k}_2 z_2^{jk} = \sup_{\vp_2 \in \Delta_2} \sum_{j,k} \pi^{*,j}_1 \pi^k_2 z_2^{jk}. $$
Hence, ${\vp}^*_i$ is a best-response for player $i$, given that the other player uses ${\vp}^*_{-i}$. While the classical theorem of Nash shows that a mixed NEP is always available, $H$ may have zero (strictly competitive case), one (standard case) or two (coordination case) pure NEP's \cite{CalvoArmengol}.

Thus, to establish existence of NEP in a multi-period game, one must consider {mixed} strategies. Furthermore, since multiple equilibria are possible, an \emph{equilibrium selection} rule is needed. In the context of the $CO_2$ emission game,  because the $(P,X$)-prices are stochastic, it is impossible to \emph{a priori} rule out some of the above scenarios for all possible state variable realizations.
In particular, the case of  the anti-coordination ``battle-of-the-sexes'' or ``chicken'' game is likely to appear when the electricity-carbon spread is slightly positive. In such a situation, each of the players will have an incentive to emit; however, if the price impact is strong enough, it is not profitable for \emph{both} of them to consume permits. As a result, two pure Nash equilibria are possible whereby one producer {yields} the market to the other.

The communication device $\gamma$, introduced in one-shot games by \cite{Aumann87,Myerson86}, provides a general method for describing such coordination while maintaining the non-cooperative game setting.  
\begin{deff}
A Markovian correlated equilibrium point for the switching game is a Markov communication device $\gamma : (s,p,x,\vze) \to \Delta_4$ inducing admissible stage strategy profiles $\vxi^*(t; \gamma) = (\mu_1(t, \gamma), \mu_2(t, \gamma))$ such that $\forall (s, p_0, x_0, \vze)$
(recall definition of $V_i$ in \eqref{eq:expected-profit})
\begin{align}\label{eq:production-ce}
\left\{ \begin{aligned}
V_1(s, p_0, x_0, \vze; u^*_1, u^*_2) &\ge V_1(s, p_0, x_0, \vze; u_1, u^*_2) \qquad \forall u_1\in \U_1, \\
V_2(s, p_0, x_0, \vze; u^*_1,u^*_2) &\ge V_2(s, p_0, x_0, \vze; u^*_1, u_2) \qquad \forall u_2 \in \U_2.
\end{aligned}\right.
\end{align}
The resulting game values are denoted as $V_i(s, p_0, x_0, \vze; \gamma)$.
\end{deff}

The meaning of the {correlated equilibrium} in \eqref{eq:production-ce} is that conditional on the private signal sequence, neither player has an incentive to deviate from the prescribed action. Therefore, given $\mu_i(t, \gamma)$ and market information $\F_t$, it is optimal to take $u_i(t) = \mu_i(t, \gamma)$.
Note that in \eqref{eq:production-ce}, even if a player chooses to deviate from the recommendation $\mu_i(t, \gamma)$ she continues to receive future signals $\mu_i(s, \gamma)$, $s>t$ and therefore information about the implied strategy of the other player. Existence of CEP of switching games will be established in Theorem \ref{thm:iter-game}. We will also
provide a recursive construction of $V_i(t, \cdot)$ in terms of conditional expectations of $V_i(t+1,\cdot)$ and one-shot $2 \times 2$ games. This allows for a solution method, detailed in Section \ref{sec:numerical}, analogous to the dynamic programming paradigm for ordinary stochastic control problems. Finally, we will show that CEP of switching game induces rational behavior at each stage, i.e.\ matches with a CEP of a one-stage sub-game.

\begin{rem}
A related concept of competitive equilibrium in the industrial organization literature is that of a stage Stackelberg game \cite{Basar}. In a Stackelberg game, at each stage one player is the {leader} and has priority in making decisions; the second player then follows. This description corresponds closely to the preferential mechanism of equilibrium selection which always favors the leader.
\end{rem}

Economically, the third party implementing the correlation could be a government regulator,  market watchdog, or just a proxy for market frictions that make one equilibrium most preferable. Thus, no inherent collusion is required and the game is still non-cooperative. If a regulator is involved, a socially beneficial correlation can be selected. For instance, a ``utilitarian'' communication device maximizes the (weighted) sum of the firms continuation values so that the producers as a whole have best economic health.  Alternatively, a ``green'' device minimizes $C O_2$ emissions. Finally, a ``preferential'' communication mechanism can endogenously emerge without a third party due to extra advantages available to a given player (e.g.\ due to preferential regulatory treatment or other externalities).

%

\begin{rem}
A variety of correlated decision-making is possible in sequential games \cite{Myerson86}. Here we focus on the stepwise weak communication device whereby the players and the regulator communicate before each stage; such a formulation allows the most flexibility and fits our economic description. However, in practice much weaker correlation could suffice. For instance, players can agree at date 0 to use the preferential-$i$ correlation law which means that in any ``tie-break'' case, player $i$ ``wins''. Once this rule is fixed, no further communication would be necessary. Similarly, if $\gamma$ is such that the implied strategy ${\vp}_{-i}(t,\gamma)|\mu_i(t,\gamma)$ of the other player is always pure, then a public randomization is sufficient at each step and no private signals are needed. Any mixture of NEP's is also a CEP and therefore except for the strictly-competitive games, one may always find correlation devices that correspond to pure Nash equilibria, obviating the need for randomization (either by players or regulator).
 
\end{rem}

\section{Sequential Stopping Game}\label{sec:seq-games}


%

Our analysis of the switching game will consist of building up the solution in several steps. We begin with analyzing the single-agent objective.
Next, in Section \ref{sec:stopping-game} we move on to the one-shot non-zero-sum stopping game that is built iteratively from the one-period $2 \times 2$  games, following the methods of \cite{RamseySzajowski08}. Finally, in Section \ref{sec:multiple-stopping-game} we describe the sequential stopping game that in the limit is shown in Section \ref{sec:seq-games} to coincide with our original model in \eqref{eq:expected-profit}.

\subsection{Single Producer Problem}\label{sec:single-player}
Before tackling the stochastic duopoly game, let us briefly review the solution of the single-player model. Since the control $u_i(t)$ takes on a finite number of values, we have an optimal switching model that can be viewed as a sequence of optimal stopping problems. Such models (including price impact) were studied in \cite{CarmonaLudkovskiSwitching,LudkovskiFlex}.

Let us consider the optimization for producer 1. For the remainder of this section we fix a production schedule $u_2$ of the second producer, as well as a communication device $\gamma$ that sends private signals $\mu_1(t,\gamma)$ to player 1.
In the single-producer problem, the objective is to maximize the expected profit
\begin{align}\label{eq:single-agent-obj}
\sup_{(u(t)) \in \U_1(\gamma)} \E^{(u,u_2,\gamma)} \left[ \sum_{t=0}^{T-1} \left\{(a_1 P_t - b_1 X_t - c_1) u(t) \, \dt - K_{\{1, u(t-), u(t)\}} \right\} \right].
\end{align}

Consider initial conditions $P_s = p, X_s = x, u_2(s) = \zeta_2$ and let $V(s, p,x,\zeta_2)$ be the value function corresponding to \eqref{eq:single-agent-obj} conditional on starting in the ``on''-production regime, and $W(s,p,x,\zeta_2)$ the value function starting offline.  Furthermore, using same initial conditions define recursively \begin{align}\label{eq:1-player-recursion}
\left\{\begin{aligned}
V^0(s,p,x,\zeta_2) & = \E^{(1,u_2,\gamma)}\left[ \sum_{t=s}^{T-1} (a_1 P_t - b_1 X_t - c_1) \, \right], \qquad\quad\text{as well as}\qquad W^0(s, p,x,\zeta_2) = 0; \\
V^n(s,p,x,\zeta_2) & = \sup_{\tau \in \sD^1(s)} \E^{(1,u_2,\gamma)}\left[ \sum_{t=s}^{\tau-1} (a_1 P_t - b_1 X_t - c_1) + (W^{n-1}(\tau, P_\tau, X_\tau, u_2(\tau)) - K_{\{1,1,0\}}) \right]; \\
W^n(s,p,x,\zeta_2) & = \sup_{\tau \in \sD^1(s)} \E^{(0,u_2,\gamma)}\left[ V^{n-1}(
\tau, P_\tau, X_\tau,u_2(\tau)) - K_{\{1,0,1\}} \right], \qquad n \ge 1\end{aligned}\right.
\end{align}
where under $\PP^{(i,u_2,\gamma)}$ the drift of the carbon allowance price is $f(i, u_2(t) )$.

\begin{prop}\label{prop:single-agent}
Let $\mathcal{U}^n_1 \triangleq \{ u \in \mathcal{U}_1 : u \text{ has at most n switches} \}$. Then, $$V^n(s,p,x, \zeta_2) = \sup_{(u(t)) \in \mathcal{U}^n_1, u(s-) = 1} \E^{(u,u_2,\gamma)}\left[ \sum_{t=s}^{T-1} \left\{(a_1 P_t - b_1
X_t - c_1) u(t) \, \dt - K_{\{1, u(t-), u(t)\}}  \right\} \right],$$ and as $n\to \infty$, $V^n(s,p,x, \zeta_2) \to V(s,p,x, \zeta_2)$, $W^n(s,p,x, \zeta_2) \to W(s,p,x, \zeta_2)$ uniformly on compacts.
\end{prop}

\begin{proof}
This is an analogue of \cite[Theorem 1]{CarmonaLudkovskiSwitching}. Compared to our earlier work, the only new feature is that the payoffs to producer 1 are randomized. Indeed, from her perspective, the strategy of player 2 (implied through the private signal $\mu_1(t,\gamma)$) may be mixed. Consequently, her continuation value is unknown at decision time, depending as it is on the action of player 2. Formally, allowing for a relaxed switching control $p^1_s$ at date $s$ (representing probability of being on) the dynamic programming principle implies that in \eqref{eq:1-player-recursion}
%
\begin{align*}
V^n &(s,p,x, \zeta_2) 
   =   \E^{\gamma(s)} \Biggl[ \sup_{p^1_s \in [0,1]} \Bigl\{ p^1_s(a_1 p - b_1 x - c_1) - (1-p^1_s) K_{\{1,1,0\}} \\
&+  \E^{\mu_1(s,\gamma)}\Bigl[ p^1_s p^2_s V^{n}(s+1, P_{s+1}, X^{(0,1)}_{s+1}, 1) + p^1_s(1-p^2_s)V^{n}(s+1, P_{s+1}, X^{(0,0)}_{s+1}, 0) \\
&+  (1-p^1_s)p^2_s W^{n-1}(s+1, P_{s+1}, X^{(0,1)}_{s+1}, 1) + (1-p^1_s)(1-p^2_s) W^{n-1}(s+1, P_{s+1}, X^{(0,0)}_{s+1}, 0) \Bigr] \Bigr\} \Biggr].
\end{align*}
The outer expectation is averaging over the signal $\mu_1$ whose law is specified by the communication device $\gamma$; however the decision-maker has access to $\mu_1(t, \gamma)$ and therefore makes the switching decision $p^1_s$ based on the conditional strategy $(p^2_s)|\mu_1(s, \gamma)$ of player 2. The inner optimization is linear in $p^1_s$ and therefore the optimizer must be an endpoint of $[0,1]$. Thus, as expected, given the signal we can work with pure controls, $u_1(t) \in \F_t \vee \sigma(\mu_1(t, \gamma))$. Note that from the perspective of an observer who has access only to $\F$, the strategy of both players appears randomized.


The rest of the proof proceeds exactly as in \cite{CarmonaLudkovskiSwitching} by iterating over the control decisions of producer 1 using the strong Markov property of $(P,X)$ and the Snell envelope characterization of optimal stopping problems.
\end{proof}



Proposition \ref{prop:single-agent} shows that the solution to \eqref{eq:single-agent-obj} can be represented in terms of the sequence $(V^n,W^n)$ which correspond to optimal stopping problems defined in \eqref{eq:1-player-recursion}. Taking the limit $n \to \infty$ we obtain 
\begin{cor}
$(V,W)$ satisfy the coupled dynamic programming equation:
\begin{align*}
\left\{ \begin{aligned}
V(s,p,x,\zeta_2) & = \sup_{\tau \in \sD^1(s)} \E^{(1,u_2,\gamma)}\left[ \sum_{t=s}^{\tau-1} (a_1 P_t - b_1 X_t - c_1) \, \dt + (W(\tau, P_\tau, X_\tau, u_2(\tau)) - K_{\{1,1,0\}}) \right], \\
W(s,p,x,\zeta_2) & = \sup_{\tau \in \sD^1(s)} \E^{(0,u_2,\gamma)}\Bigl[ V(\tau, P_\tau, X_\tau, u_2(\tau)) - K_{\{1,0,1\}} \Bigr]. 
\end{aligned} \right.
\end{align*}
Moreover, an optimal strategy $u^*_1 \in \U_1$ exists.
\end{cor}

\subsection{Correlated Equilibria in Non-Zero-Sum Stopping Games}\label{sec:stopping-game}
In this section we recall existing results on two-player non-zero sum stopping games in {discrete time and finite horizon}. Let $\Z \equiv (Z^{jk}_i(t))$, $i\in \{1,2\}, j,k \in \{0,1\}$ be a octuple of bounded $(\F_t)$-adapted stochastic processes.
Player $i \in \{1,2\}$ optimizes the reward
\begin{align}\label{eq:stopping-game-reward}
\tJ_i(s,\tau_1, \tau_2) \triangleq\left(\sum_{t=s}^{\tau_i \wedge \tau_{-i}-1} Z_i^{00}(t) \right) + Z_i^{10}(\tau_i) 1_{\{ \tau_i < \tau_{-i}\}} + Z_i^{01}(\tau_{-i}) 1_{\{ \tau_{-i} < \tau_i \}} + Z^{11}_i(\tau_i) 1_{\{ \tau_i = \tau_{-i} \}},
\end{align}
by choosing the (randomized) $(\mathcal{F}_t)$-stopping time $\tau_i \le T$. In words, $Z^{00}_i$ is the ongoing reward for staying in the game, $Z^{10}_i$ is the reward if the player stops first; $Z^{01}_i$ is the reward if the other player stops first and $Z^{11}_i$ is the reward if both players stop simultaneously. Thus, continuing is associated with action `0' and stopping with action `1'.

The Dynkin zero-sum stopping game corresponds to
$Z^{10}_1 = -Z^{01}_2, Z^{01}_1 = -Z^{10}_2, Z^{11}_1= -Z^{11}_2$ and was recently fully analyzed  by \cite{EkstromPeskir08}. Also, the monotone cases $Z^{01}_i \le Z^{11}_i \le Z^{10}_i$ $\bar{\PP}$-a.s.\ (where both players prefer to stop late) and $Z^{01}_i \ge Z^{11}_i \ge Z^{10}_i$ were considered by Ohtsubo  \cite{Ohtsubo87}. In these special cases, a unique pure Markov NEP exists. The fundamental result of \cite{Ohtsubo87,Ohtsubo91} characterizes game value functions $(V_1, V_2)$ for $\Z$ as  a pair of $\F$-adapted processes satisfying $\E[ \sup_{0 \le t \le T} V_i(t) ] < \infty$, $V_i(T) = Z^{11}_i(T)$ and for all $0 \le t \le T$
\begin{align}\label{eq:val-H-recursion}
(V_1(t), V_2(t)) \in \mathcal{E} \begin{pmatrix} ( \E[ V_1(t+1) | \F_t]+Z^{00}_1(t) , \E[ V_2(t+1) | \F_t]+Z^{00}_2(t) )   & (Z^{01}_1(t), Z^{01}_2(t)) \\ (Z^{10}_1(t), Z^{10}_2(t)) & (Z^{11}_1(t), Z^{11}_2(t)) \end{pmatrix},
\end{align}
%
%
where $\mathcal{E}(H)$ is the set of game values corresponding to NEPs of $H$. 
This reduces computation of game values to iterative solution of one-shot $2 \times 2$ games, in complete analogy to standard dynamic programming. We seek a similar result for the switching game, see \eqref{eq:val-V-recursion} below.

Without any assumptions on the structure of $\Z$ appearing in \eqref{eq:stopping-game-reward}, the existence of a pure NEP is not guaranteed. However, as shown by \cite{Ferenstein07} (see also \cite{ShmayaSolan} and references therein) a two-person stopping game always admits a mixed NEP.
Again, there is no uniqueness and we might need equilibrium selection.
%
Let $\gamma$ be an $(\mathcal{F}_t)$-adapted stochastic process taking values in
$\Delta_4$. Define the dependent randomized stopping rules (cf.~\eqref{def:signals})
$$\left\{ \begin{aligned} \tau_1(\gamma) & \triangleq \inf\{ t : \eta'(t) \le \gamma_{10}(t)+\gamma_{11}(t) \}, \\ \tau_2(\gamma) & \triangleq \inf\{ t : \eta'(t) \le \gamma_{01}(t) + \gamma_{11}(t) \},  \end{aligned}\right.\qquad \eta'(t) \sim Unif[0,1] \quad\text{   i.i.d.}.$$
Thus, conditional on the game still continuing, the stage-$t$ payoff to player $i$ is $\sum_{j,k}\gamma_{jk}(t) Z^{jk}_i(t)$ and total expected payoff is
\begin{align}
\E^\gamma \left[\tJ_i(s,\tau_1(\gamma), \tau_2(\gamma)) \right] = \E \left[\sum_{t=s}^{T-1} \sum_{j,k} \left\{ \Bigl( \prod_{r=s}^{t-1} \gamma_{00}(r) \Bigr) \gamma_{jk}(t) Z^{jk}_i(t) \right\} \right].
\end{align}
%
%

As before, correlation is implemented through private signals $\mu_i(t,\gamma)$ and a correlated equilibrium of $\Z$ is a
%
communication device $\gamma$ inducing a stopping strategy profile $\vtau(\gamma) \triangleq (\tau_1(\gamma), \tau_2(\gamma)) \in \sD^1 \times \sD^2$ such that for $i = 1,2$ and all $0 \le  t < T$
\begin{align}\label{def:ce-stopping}
V_i(t; \gamma, \Z) \triangleq \E^{\gamma}[ \tJ_i(t, \vtau(\gamma)) | \F_t] \ge \E^{\gamma}[ \tJ_i(t, \tilde{\tau}_i, \tau_{-i}(\gamma) ) | \F_t], \quad \forall \tilde{\tau}_i \in \sD^i(t).
\end{align}
Observe that given a device $\gamma$ leading to a CEP, it must be  that
\begin{multline}\label{eq:stopping-problem-Gi}
V_i(t; \gamma, \Z) = \sup_{\tau \in \sD^i(t)} \E^{\gamma}\Bigl[ \left(\sum_{s=t}^{(\tau \wedge \tau_{-i}(\gamma))-1} Z_i^{00}(s) \right) + Z_i^{10}(\tau) 1_{\{ \tau < \tau_{-i}\}} \\ + Z_i^{01}(\tau_{-i}) 1_{\{ \tau_{-i} < \tau \}} + Z^{11}_i(\tau) 1_{\{ \tau = \tau_{-i} \}} \big| \mathcal{F}_t \Bigr]
\end{multline}
which is a standard optimal stopping problem for player $i$ in the enlarged filtration $\G^i$.

\begin{lemm}{\cite[Theorem 2.3]{RamseySzajowski08}}
Consider a CEP with communication device $\gamma$ of a stopping game $\Z$. Then for all $t \in \{0,1,\ldots, T-1\}$ we have
\begin{align}\label{eq:ce-recursion}
\left\{ \begin{aligned}
  \gamma_{00}(t) (\E[ V_1(t+1)|\F_t] + Z^{00}_1(t))  + \gamma_{01}(t) Z^{01}_1(t)  & \ge \gamma_{00}(t) Z^{10}_1(t)  + \gamma_{01}(t) Z^{11}_1(t)); \\
  \gamma_{00}(t)(\E[ V_2(t+1)|\F_t] + Z^{00}_2(t))  + \gamma_{10}(t) Z^{10}_2(t)  & \ge \gamma_{00}(t) Z^{01}_2(t)  + \gamma_{10}(t) Z^{11}_2(t) ; \\
  \gamma_{10}(t) Z^{10}_1(t)   + \gamma_{11}(t) Z^{11}_1(t)  & \ge \gamma_{10}(t)(\E[ V_1(t+1)|\F_t] + Z^{00}_1(t))  + \gamma_{11}(t) Z^{01}_1(t) ; \\
  \gamma_{01}(t) Z^{01}_2(t)  + \gamma_{11}(t)p Z^{11}_2(t)  & \ge \gamma_{01}(t)(\E[ V_2(t+1)|\F_t] + Z^{00}_2(t))  + \gamma_{11}(t) Z^{10}_2(t) . \end{aligned}
  \right.
\end{align}
\end{lemm}
Lemma 1 shows that a CEP of the stopping game is rational at each stage. For instance, the first inequality in \eqref{eq:ce-recursion} means that conditional on player 1 signal being `continue', the expected payoff to player 1 from continuing (the right-hand-side) is better than the expected payoff from stopping. In either scenario, player 2 implements the conditional strategy $\vp_2(t,\gamma) |_{ \mu_1(t,\gamma) = 0} = \bigl(\frac{\gamma_{0 0}(t)}{\gamma_{00}(t)+\gamma_{01}(t)}, \frac{\gamma_{01}(t)}{\gamma_{00}(t)+\gamma_{01}(t)} \bigr)$.

As shown by \cite[Theorem 2.4]{RamseySzajowski08}, any finite-horizon stopping game with bounded payoffs admits a CEP; in fact outside the zero-sum and monotone cases we expect that a large number of CEPs are possible.
It is convenient to think of communication device $\gamma$ leading to a CEP in \eqref{def:ce-stopping} as a measurable selector of local correlated equilibrium points in the one-shot $2 \times 2$ games. Thus, let $\Gamma : \bT \times \Omega \times \R^{2 \times 2 \times 2} \to \Delta_4$ be a measurable map such that for any $2 \times 2$ game $H$, $\Gamma(t, \omega, H)$ is a CEP of $H$. Then using $\Gamma$, one may construct a communication device $\gamma$ by inductively using the CEP $\Gamma(t, \omega, H(t, \omega))$, where $H(t,\omega)$ is the right-hand-side in \eqref{eq:val-H-recursion}, and proceeding back in time. Observe that for most $H$'s, $\Gamma(\cdot, H)$ is simply the unique NEP available, so that the selection feature is ``silent'', and the device is only really activated when considering the coordination game. With this perspective in mind, we call a \emph{correlation law} $\Gamma$ a communication device which is based on the same local criterion (for instance ``minimize today's emissions'' or ``maximize today's value of player 1'').

\subsection{Recursive Construction}\label{sec:multiple-stopping-game}
We return to the emissions market duopoly setup. The emission schedules of the two agents are interpreted as a sequence of regime-changes. Thus, the single-stopping game in the previous section is viewed as the sub-game for making the next regime-switch. The stopping game in Section \ref{sec:stopping-game} is accordingly denoted as a $(1,1)$-fold switching game and we now will consider $(n,m)$-fold switching games with game value functions $V^{n,m}$. These games have a restricted set of possible production strategies; namely the total number of regime switches over the game horizon is bounded by $n$ and $m$, respectively. Using the Markov property of the game state and actions it is not surprising that these various switching games are related to each other.

In terms of the notation of Section \ref{sec:stopping-game}, we identify the running profit with $Z^{00}_i(t) = (a_i P_t - b_i X_t - c_i)u_i(t)$ and the other $Z_i^{jk}$'s with various game continuation-values.
%
For the remainder of the section, we make a standing assumption that a communication device $\gamma$ is chosen and fixed.
Let us fix an initial state $P_s = p, X_s = x$ and initial production regime $\vze = (\zeta_1, \zeta_2)$. Define a double cascade of stopping games indexed by $n$ and $m$ via
\begin{align}\label{eq:hJ}
\hJ_i^{n,m}(s,p,x,\vze) & \triangleq V_i(s; \gamma, \tilde{\Z}^{n,m}(\vze)), \qquad n,m \ge 1
\end{align}
which uses the notation of \eqref{def:ce-stopping} based on the recursive payoff structure
\begin{align}\label{eq:tildeZ-def}\left\{ \begin{aligned}
(\tilde{Z}^{n,m})^{00}_i(t, \vze) & = (a_i P_t - b_i X_t - c_i)\zeta_i; \\
(\tilde{Z}^{n,m})^{01}_i(t, \vze) & = \hJ_i^{n,m-1}(t, P_t, X_t,\zeta_1,1-{\zeta}_2) -1_{\{i=2\}}K_{\{2,\zeta_2, 1-{\zeta_2}\}};\\
(\tilde{Z}^{n,m})^{10}_i(t, \vze) & = \hJ_i^{n-1,m}(t, P_t, X_t,1-{\zeta_1},\zeta_2) - 1_{\{i=1\}}K_{\{1,\zeta_1, 1-{\zeta_1}\}};\\
(\tilde{Z}^{n,m})^{11}_i(t, \vze) & = \hJ_i^{n-1,m-1}(t, P_t, X_t,1-{\zeta_1}, 1-{\zeta_2}) - K_{\{i,\zeta_i, 1-{\zeta_i}\}}.\end{aligned}
\right. \end{align}
The boundary cases are first
%
%
%
\begin{align*}
& \hJ_i^{0,0}(s, p,x,\vze) \triangleq \E^{\vze} \left[\sum_{t=s}^{T-1} (a_i P_t - b_i X_t - c_i) \zeta_i \right];
\end{align*}
next,  $\hJ_1^{n,0}(s,p,x,\vze)$  and $\hJ_2^{0,m}(s,p,x,\vze)$ are identified with the single-player optimization problems as in \eqref{eq:1-player-recursion} (keeping the emission regime of the other player fixed at $\zeta_{-i}$). Finally, we take
$$
\hJ_2^{n,0}(s,p,x,\vze) = \E^{(u^{n,*}_1,\zeta_2,\gamma)}\left[ \sum_{t=s}^{T-1} (a_2 P_t - b_2 X_t - c_2) \zeta_2 \right]
$$
where $u^{n,*}_1$ is an optimal control for the problem defining $\hJ_1^{n,0}$,  and similarly for $\hJ_1^{0,m}(s,p,x,\vze)$.

%

\subsection{Switching Game Equilibrium as Sequential Stopping Game Equilibrium}
We now proceed to \emph{glue} the sequential stopping games of $\hJ_i^{n,m}$ and re-interpret the latter as value functions of a switching game. For $n \ge 0$ denote by $\U^{n}_i \subset \U_i$ the set of all production strategies for player $i$ with at most $n$ switches. Consider the restricted repeated game with payoffs \eqref{eq:expected-profit} where we require $u_1 \in \U_1^{n}$ and $u_2 \in \U_2^{m}$, so that the first producer may change her production regime at most $n$ times, and the second producer at most $m$ times.

Our first task is to obtain a switching-game CEP that matches the definition of $\hJ^{n,m}$. To do so we pick a correlation law $\Gamma$; $\Gamma$ gives rise to a CEP of any stopping game, in particular it leads to well-defined game values $\hJ_i^{n,m}$ in \eqref{eq:hJ}. We now construct a communication device $\gamma^{n,m}$ for the $(n,m)$-switching game. Let $k_i(t)$ be the number of production switches used by player $i$ by stage $t$. The device $\gamma^{n,m}(t)$ at stage $t$ is taken to be $\Gamma\left(t, \omega, \tilde{\Z}^{n-k_1(t),m-k_2(t)}(\vxi(t)) \right)$ defined in terms of \eqref{eq:tildeZ-def} and the latest regime $\vxi(t)$. Note that the overall $\gamma^{n,m}$ is no longer Markovian since it has memory of the number of switches made by each player, which is necessary in the constrained game. The above construction is well-defined for all paths of $(P,X, \vxi)$, even outside equilibrium.

Using $\gamma^{n,m}$ we proceed to construct switching controls $u_i^{n,m}$ for the $(n,m)$-switching game. To simplify notation we write $\utau^{n,m} = \tau_1(\gamma^{n,m}) \wedge \tau_2(\gamma^{n,m})$ which is interpreted as the equilibrium first stopping time for the game defined by \eqref{eq:hJ} under the correlation law $\Gamma$. Given the starting production regime $\vze = (\zeta_1, \zeta_2)$, let us define the switching controls $u_i^{n,m}(s)$ for this game by
\begin{align}\notag
u_1^{n,m}(s) &= \zeta_1 \qquad\text{for}\quad s < \utau^{n,m}; \\ \label{eq:opt-control}
u_1^{n,m}(s) &= \left\{ \begin{aligned} 1-{\zeta_1} & \text{ for}\quad \utau^{n,m} \le s < \utau^{n-1,m} & \quad\text{when } \tau^{n,m}_1 < \tau^{n,m}_2; \\
\zeta_1 & \text{ for}\quad \utau^{n,m} \le s < \utau^{n,m-1} & \quad\text{when } \tau^{n,m}_2 < \tau^{n,m}_1; \\
1-{\zeta_1} & \text{ for}\quad \utau^{n,m} \le s < \utau^{n-1,m-1} & \quad\text{when } \tau^{n,m}_1 = \tau^{n,m}_2, \end{aligned}\right. \\ \notag ... \text{ and so on,}&
\end{align}
and   similarly for  $u_2^{n,m}(t)$.
In words, $u_i^{n,m}$ keeps track of the production regime of the $i$-th agent following the decision rules defined sequentially by descending through the family of the $\hJ^{n,m}$\nobreakdash-\hspace{0pt}stopping subgames (one stopping game at a time).  Then by definition of \eqref{eq:hJ} we have $u_1^{n,m} \in \U_1^n$ and $u_2^{n,m} \in \U_{2}^m$. It can also be seen through an easy induction argument that
\begin{align}\label{eq:Vnm-two-ways}
\hJ^{n,m}_i(s,p,x,\vze) = V_i( s,p,x,\vze; \vxi^{n,m}),
\end{align}
so that the switching control $\vxi^{n,m}$ of \eqref{eq:opt-control} allows to achieve the game values $\hJ^{n,m}$ defined recursively in \eqref{eq:hJ}.
Moreover, the next theorem shows
 that the pair $(u_1^{n,m}, u_2^{n,m})$ is in fact a correlated equilibrium (using correlation device $\gamma^{n,m}$) for the game \eqref{eq:expected-profit} over the control set $\mathcal{U}_1^{n} \times \mathcal{U}_2^{m}$.

\begin{thm}\label{thm:iter-ce}
For all $n > 0$ and $u_1 \in \U_1^n$ we have $V_1(t, \cdot; u_1^{n,m}, u_2^{n,m}) \ge V_1(t, \cdot; u_1, u_2^{n,m})$. Similarly for all $m>0$ and $u_2 \in \U_2^m$ we have $V_2(t, \cdot; u_1^{n,m}, u_2^{n,m}) \ge V_2(t, \cdot; u_1^{n,m},u_2)$.
\end{thm}
\begin{proof}
The idea of the proof is to make use of the Markov structure of our problem and apply induction. The other key tool is that given $\gamma^{n,m}$, we can look at one player at a time which essentially reduces to a single-player problem studied before, see \eqref{eq:stopping-problem-Gi}.

Due to symmetry, it suffices to prove the result for player 1. When $m=0$ the other player cannot act, the game becomes trivial and Theorem \ref{thm:iter-ce} is just a re-statement of Proposition \ref{prop:single-agent}. Conversely, when $n=0$, the first player cannot act and there is nothing to prove. Using induction we assume that the theorem has been shown for the pairs $(n-1,m-1)$, $(n,m-1)$ and $(n-1,m)$; let us show it for the case $(n,m)$. Given an arbitrary $u_1 \in \U_1^{n}$, write it as $u_1 = (\tau^1, \hat{u}_1)$ where $\hat{u}_1 \in \U_1^{n-1}$ denotes the remainder of $u_1$ after the first switch time $\tau^1$.  Let $\tau^{2,*} \equiv \tau^2(\gamma^{n,m})$ be the first switch for the second player dictated through $\gamma^{n,m}$. Define $\utau = \tau^1 \wedge \tau^{2,*}$. Also for notational convenience we omit all the arguments of $V^{n,m}$ except for the time variable.
Then the strong Markov property of $(P,X)$ and the way $u_2^{n,m}(t)$ was constructed show that
\begin{multline*}
\E^{(u_1, u_2^{n,m},\gamma^{n,m})}\left[ \sum_{s=\utau}^{T-1} (a_1 P_s - b_1 X_s - c_1) \hat{u}_1(s)
\right] = \E^{(u_1, u_2^{n,m},\gamma^{n,m})}\Bigl[ V_1( \utau; \hat{u}_1, u_2^{n-1,m}) 1_{\{ \tau^1 < \tau^{2,*}\}} \\ + V_1( \utau; \hat{u}_1, u_2^{n,m-1}) 1_{\{ \tau^1 > \tau^{2,*}\}} + V_1( \utau; \hat{u}_1, u_2^{n-1,m-1}) 1_{\{ \tau^1 = \tau^{2,*}\}} \Bigr].
\end{multline*}
Conditioning on $\tau^1$ and $\tau^{2,*}$ we therefore have $V_1(t ; u_1, u_2^{n,m}) =$
\begin{align*}
 & \E^{(u_1, u_2^{n,m},\gamma^{n,m})}\Bigl[ \left( \sum_{s=t}^{\utau-1} (a_1 P_s - b_1 X_s - c_1) u_1(t) \right) + \left(\sum_{s=\tau^1}^{T-1} (a_1 P_s - b_1 X_s - c_1) \hat{u}_1(s) \right) 1_{\{\tau^1 < \tau^{2,*}\}} \\
 & + \Bigl( \sum_{s=\tau^{2,*}}^{T-1}  (a_1 P_s - b_1 X_s - c_1) \hat{u}_1(s) \Bigr)1_{\{\tau^1 > \tau^{2,*}\}}  + \left( \sum_{s=\tau^1}^{T-1} (a_1 P_s - b_1 X_s - c_1) \hat{u}_1(s) \right)1_{\{\tau^1 = \tau^{2,*}\}} \Bigr] \\
 = &\E^{(u_1, u_2^{n,m},\gamma^{n,m})}\Bigl[ \left( \sum_{s=t}^{\utau-1} (a_1 P_s - b_1 X_s - c_1) u_1(t) \right)+ V_1( \utau; \hat{u}_1, u_2^{n-1,m}) 1_{\{ \tau^1 < \tau^{2,*}\}}  \\ & \qquad + V_1( \utau; \hat{u}_1, u_2^{n,m-1}) 1_{\{ \tau^1 > \tau^{2,*}\}} + V_1( \utau; \hat{u}_1, u_2^{n-1,m-1}) 1_{\{ \tau^1 = \tau^{2,*}\}} \Bigr]\\
\intertext{by induction hypothesis we have the inequality}
 \le &\E^{(u_1, u_2^{n,m},\gamma^{n,m})}\Bigl[ \left( \sum_{s=t}^{\utau-1} (a_1 P_s - b_1 X_s - c_1) u_1(t)  \right) + V_1(\tau^1; {u}_1^{n-1,m}, u_2^{n-1,m}) 1_{\{\tau^1 < \tau^{2,*}\}} \\
& \qquad + V_1(\tau^{2,*}; {u}_1^{n,m-1}, u_2^{n,m-1}) 1_{\{\tau^1 > \tau^{2,*}\}}  +  V_1(\tau^1; {u}_1^{n-1,m-1}, u_2^{n-1,m-1}) 1_{\{\tau^1 = \tau^{2,*}\}} \Bigr] \\
 \le &\sup_{\tau^1  \in \sD^1(t)} \E^{(u_1(t), u_2^{n,m},\gamma^{n,m})}\Bigl[ \left( \sum_{s=t}^{\utau-1} (a_1 P_s - b_1 X_s - c_1) u_1(t)  \right) + V_1(\tau^1; {\vxi}^{n-1,m} ) 1_{\{\tau^1 < \tau^{2,*}\} } \\
& \qquad + V_1(\tau^{2,*}; {\vxi}^{n,m-1}) 1_{\{\tau^1 > \tau^{2,*}\}}  +  V_1(\tau^1; {\vxi}^{n-1,m-1}) 1_{\{\tau^1 = \tau^{2,*}\}} \Bigr] \\
 = &V_1(t; \vxi_1^{n,m}),
\end{align*}
where the last line uses the relationship \eqref{eq:Vnm-two-ways}, the construction of the stopping game defining $V^{n,m}$ in \eqref{eq:tildeZ-def}, and property \eqref{eq:stopping-problem-Gi}.
\end{proof}


The above construction leads to the key result of this section that characterizes CEP of switching games, establishes their existence, and gives a recursive formula for the resulting game value functions. For a $2 \times 2$ game $H$ defined in \eqref{eq:H-def} and correlation device $\gamma$ we denote the respective game values as

\begin{thm}\label{thm:iter-game}
Fix a  correlation law $\Gamma$. Then $\Gamma$ gives rise to a CEP of the switching game \eqref{eq:expected-profit}. Moreover, the corresponding value functions $V_i(t,P_t,X_t, \vze; \Gamma)$ solve
\begin{align}\label{eq:val-V-recursion}\left\{ \begin{aligned}
V_1(t,P_t, X_t,\vze) & = \gamma_{00}(t)Y_1(\zeta_1, \zeta_2) + \gamma_{01}(t)Y_1(\zeta_1, 1-{\zeta_2}) \\ & \qquad + \gamma_{10}(t)(Y_1(1-{\zeta_1}, {\zeta_2})-K_{1,\zeta_1,1-{\zeta}_1}) + \gamma_{11}(t) (Y_1(1-{\zeta}_1, 1-{\zeta_2})-K_{1,\zeta_1,1-{\zeta}_1})\\
V_2(t,P_t, X_t,\vze) & = \gamma_{00}(t)Y_2(\zeta_1, \zeta_2) + \gamma_{01}(t)(Y_2(\zeta_1, 1-{\zeta_2})-K_{2,\zeta_2, 1-{\zeta}_2}) \\
& \qquad + \gamma_{10}(t)Y_2(1-{\zeta_1},\zeta_2) + \gamma_{11}(t)(Y_2(1-{\zeta_1}, 1-{\zeta_2})-K_{2,\zeta_2, 1-{\zeta}_2})
\end{aligned}\right.\end{align}
where $Y_i(t,\vze) \triangleq \E^{\vze}[ V_i(t+1, \vze) | \F_t] + (a_i P_t - b_i X_t - c_i)\zeta_i$ and we have omitted the dependence on $t$.
The equilibrium controls can be taken as $u^*_i \equiv u^{T,T}_i$, as defined in \eqref{eq:opt-control}.
\end{thm}

\begin{proof}
We wish to take $n,m \to \infty$ in Theorem \ref{thm:iter-ce}. Because for $n>m$, $\U^m \subseteq \U^n$, it follows that for a fixed $m$, $\hJ_1^{n,m}$ is increasing in $n$ (and for a fixed $n$, $\hJ_2^{n,m}$ is increasing in $m$). For our discrete-time game, at most $T$ regime switches are possible for each player. Therefore $u_i^* \in \U_i^T$ and it follows that  $\hJ^{n,n}_i \equiv V_i$ for all $n >T$.
In particular, a switching CEP based on $\Gamma$ results by using $\gamma^{T,T}$.

 Moreover, at equilibrium at most one switch is made at any given stage due to the triangle condition on $K_{i,j,k}$. Therefore, if it is optimal to switch at stage $t$ from $\vze$ to $\vxi$, then already starting at regime $\vxi$ at $t$ and same state variables it is optimal to make no changes, so that
$V_i(t,\vxi) = \E^{\vxi}[ V_i(t+1,\vxi) | \F_t] + (a_i P_t - b_i X_t - c_i)u_i$ for that scenario. Combining these facts with the form of \eqref{eq:val-H-recursion} and dropping the constraints on the number of switches, we may express all payoffs in terms of next-stage game values. The recursion \eqref{eq:val-V-recursion} is now obtained by making this substitution in \eqref{eq:val-H-recursion}.
\end{proof}


%



\section{Numerical Implementation}\label{sec:numerical}


Theorem \ref{thm:iter-game} shows that a game value and equilibrium strategy profile can be obtained recursively by solving the 1-period 2-by-2 games in \eqref{eq:val-V-recursion}. The payoffs of those games are given in terms of \emph{conditional expectations} of next-stage game values. Therefore, a numerical implementation hinges on accurate evaluation of these expectations. Since our state-space in $(P,X)$ is continuous, it is impossible to make this computation exactly. Below we present two possible approximation approaches.

\subsection{Markov Chain Approximation Algorithm}
Our model would be simplified if the continuous state space of $(P,X)$ is discretized. Let $(\tilde{P}, \tilde{X})$ be an approximating discrete-state process with $(\tilde{P}_t, \tilde{X}_t)$ living on a finite subset $D_t \subset \R_+^2$. If the pair $(\tilde{P}, \tilde{X})$ is furthermore chosen to be again Markov, this is known as the Markov Chain Approximation (MCA) method of \cite{KushnerDupuis}. With such $(\tilde{P}, \tilde{X})$, a conditional expectation $\E[ f({P}_{t+1}, {X}_{t+1}) | {P}_t=p, {X}_t=x] \simeq \E[ f(\tilde{P}_{t+1}, \tilde{X}_{t+1}) | \tilde{P}_t = p, \tilde{X}_t = x]$ for any measurable function $f$ is just a weighted sum based on the transition probability matrix of $(\tilde{P}, \tilde{X})$. The backward recursion in \eqref{eq:val-V-recursion} for $\tilde{V}_i$, the corresponding approximation of $V_i$, can now be implemented directly for each stage $t$ and each possible state of $(\tilde{P}_t, \tilde{X}_t) \in D_t$. A well-known procedure constructs $(\tilde{P}, \tilde{X})$ by taking $D_t$ to be a 2-dimensional regular grid or lattice and allowing state transitions only between neighboring grid points.
Moreover, the transition probabilities of $(\tilde{P}, \tilde{X})$ are chosen so as to have \emph{local consistency} in  the first two moments with the 1-step transition densities of $(P,X)$; see \cite[Chapter 5]{KushnerDupuis}.

To use this approach in our model, one must take into account the price impact. Therefore, we construct four approximations $(\tilde{P}, \tilde{X}^{\vze})$ indexed by the possible joint production regimes $\vze \in \{0,1\}^2$ that induce different local dynamics of $\tilde{X}^{\vze}$, see \eqref{eq:X-dynamics}. In other words, our effective state variables are $(\tilde{P}, \tilde{X},\vze)$. For every possible combination $(t,\tilde{p}, \tilde{x}, \vze) \in \bT \times D_t \times \{0,1\}^2$ the relation \eqref{eq:val-V-recursion} is then solved through backward recursion. A generic convergence proof (as the grid spacing tends to zero) of this procedure for finite-horizon non-zero-sum stochastic games was obtained in \cite{Kushner08}. Note that in our model the controls $\vxi(t)$ are discrete and finite-valued and therefore all the compactness conditions in \cite{Kushner08}  for the control space are automatically satisfied.

\subsection{Least Squares Monte Carlo Approach}\label{sec:rmc}
Like classical dynamic programming, the MCA method above suffers from the curse of dimensionality. Indeed, the size of the approximating grid grows exponentially in the dimension of the state variables. In our basic model $(P,X)$ are two-dimensional; however realistic implementations are likely to take multi-dimensional factor models for $P$ and (possibly) $X$. Thus, it is helpful to seek a more robust algorithm.

A seminal idea due to \cite{Carriere96,Egloff05,Longstaff} is to use a cross-sectional regression combined with a Monte Carlo simulation to compute the relevant conditional expectations. The key step  is a global approximation of the maps $(t,p,x, \vze) \mapsto V_i(t,p,x, \vze)$ and equilibrium one-step strategies $(t,p,x,\vze) \mapsto \vxi(t,p,x,\vze)$ (based on a fixed correlation law $\Gamma$) via a random sample of $(P_t, X_t)$. The construction is iterative and backward in time.

Suppose that the current date is $t$ and we already know all the approximations $v_i(s,p,x,\vze) \simeq V_i(s,p,x,\vze)$ for $s>t$ and the corresponding equilibrium strategy profiles. Given a collection of initial points $(p^n_t, x^n_t)$, for $n=1,\ldots, N$, and an arbitrary starting emission regime $\vze = \vxi^n(t)$ we first simulate the future \emph{cashflows} on $[t+1,T]$ for each scenario $n$. This is done by iteratively updating $(p^n_{s+1}, x^n_{s+1})$ through an independent draw from the conditional law $\PP^{\vxi^n(s)}$ and then computing the equilibrium actions $u^n_i(s)$ of each player for $s = t+1, \ldots, T$ based on the estimated future game values $v_i(s,p^n_s,x^n_s, \vxi^n(s))$ and the chosen communication device $\Gamma$. If $\Gamma$ leads to randomized strategies, such a randomization is naturally implemented as part of this simulation. The realized pathwise cashflow $\vth^n_i(t+1, \vze)$ represents an empirical draw from $V_i(t+1,P_{t+1},X^{\vze}_{t+1}, \vze)$ conditional on $P_t = p^n_t, X_t = x^n_t$. We now perform a cross-sectional regression of $(\vth^n_i(t+1,\vze))_{n=1}^N$  against $(p^n_t, x^n_t)_{n=1}^N$ by using a collection of basis functions $B_\ell(t,p,x)$, $\ell=1,\cdots,r$. The regression yields the predicted continuation values
$$\hat{v}_i(t, p^n_t, x^n_t, \vze) \simeq \E^{\vze}
\left[ V_i(t+1,P_{t+1},X_{t+1},\vze) \big|\; P_t = p^n_t, X_t = x^n_t \right].$$
Finally, using $\hat{v}_i$ together with the current payoffs and switching costs and the correlation law $\Gamma$ we solve for the equilibrium game values $v_i(t, p^n_t, x^n_t, \vxi)$ for each production regime $\vxi$ by applying the stage-$t$ sub-game of Theorem \ref{thm:iter-game}.
  The computed game equilibrium also provides the map $(t,p^n_t, x^n_t, \vxi) \mapsto \vxi^*(t)$ for the equilibrium strategies. The regression results allow to further extend this to arbitrary initial condition $(t,p,x,\vxi)$. Working back to $t=0$, the final answer (which is a random variable depending on the Monte Carlo sample) is simply the average $V_i(0,p_0,x_0,\vze_0) \simeq \frac{1}{N} \sum_n v_i(0, p^n_0, x^n_0, \vxi_0)$.

The initial collection $(p^n_t, x^n_t)_{t=1}^T$ is obtained by simulation. Since, $X$ is affected by the price impact of $\vxi$, to simulate $(p^n_t, x^n_t)$ we need to select some anterior auxiliary strategy profile $\vxi^0$. While in theory $\vxi^0$ can be arbitrary, in practice it should be close to the equilibrium $\vxi^*$. Indeed, the collection $(v_i(t,p^n_t, x^n_t, \vze))_{n=1}^N$ is supposed to approximate $V_i(t, P_t, X^*_t, \vxi(t))$ where $X^*_t$ is the equilibrium $CO_2$ allowance price. Because $v_i$'s are computed by employing regression, the resulting approximation cannot be uniformly good on $\R_+^2$. From the point of view of accurate solutions, it needs to be good around the region of interest for $X^*_t$. Thus, we need most of the $x^n_t$'s to be in that (\emph{a priori} unknown) neighborhood. To overcome this difficulty, as the algorithm works back through time, the future paths $(p^n_s,x^n_{s})$, $s>t$ are re-computed using the now-available (approximately) equilibrium strategies $u^*(s)$. To further mitigate the problem,
we iteratively re-do the whole simulation and subsequent backward recursion a few times (in practice three iterations suffice), using the computed $\vxi^{*}$ from one iteration as the anterior $\vxi^0$ in the next one. The Appendix summarizes the above scheme in pseudo-code.

Selection of basis functions should reflect the expected shape of $(p,x) \mapsto V_i(t,p,x,\vze)$. A typical choice is to use low-degree polynomial basis functions, such as $p, p^2, x, x^2$, etc. In practice, $r=5-7$ basis functions and $N=32000-50000$ paths suffice. A large degree of customization, such as time-varying bases, constrained least-squares regression, variance reduction methods, etc., is possible to speed up the computations. The Appendix summarizes the above scheme in pseudo-code in  Algorithm \ref{alg:main}. It calls as a sub-routine Algorithm \ref{alg:forward-path} that carries out the forward simulations of $\vth^n_i$.
The cost of simulations in Algorithm \ref{alg:main} is $\mathcal{O}(N \cdot T^2)$ which consists of re-simulating $N$ paths on $[t,T]$ as $t$ goes from $T-1$ to zero (see Algorithm \ref{alg:forward-path}). The cost of doing regression against $r$ basis functions on each path and for each stage is $\mathcal{O}(N \cdot T \cdot r^3)$ and the cost of computing continuation values is $\mathcal{O}(N \cdot T^2 \cdot r)$.
The memory requirements from storing all the simulation paths  are $\mathcal{O}(N \cdot T )$.

%
%

\subsection{Numerical Examples}\label{sec:numerical-ex}
In this section we illustrate our analysis with a numerical case-study.
The selected model parameters are listed in Table \ref{table:params}. The example represents emission scheduling of two producers over one calendar year; all the parameters of $(P,X)$ are in annualized units and we use bi-weekly periods $T' = 26$  to model the scheduling flexibility. Note that the electricity price $P_t$ is more volatile than the $CO_2$ allowance price $X_t$; also the mean-reversion parameter $\kappa_X$ of $X$ is quite large, implying a significant price impact. In \eqref{eq:X-dynamics} we take $f( \zeta_1, \zeta_2) = \log(12 + 8 \zeta_1 + 4 \zeta_2)$, so that the mean-reversion level of $\log X$ is linear in the production regimes of producers 1 and 2, with producer 1 having more influence due to emitting twice as much carbon, $b_1 = 2b_2 \Rightarrow g_1 = 2g_2$. The stylized production/emission parameters represent a dirty ``coal'' producer 1 who has low input costs but needs lots of allowances, and a clean ``natural gas'' producer 2 who has high fixed costs but small sensitivity to allowance prices (and can generate twice as much electricity). Observe that if both producers emit simultaneously for a long period of time, then we expect $P_t \sim \bar{P} = 45, X_t \sim f(1,1) = 24$ meaning that everyone will be losing money. Therefore, extended joint emissions are not \emph{sustainable}.

\begin{table}
\begin{tabular*}{\textwidth}{lr}
\begin{minipage}{3in}
$$ \begin{array}{|cc|cc|} \hline
\kappa_X & 3 & \sigma_X & 0.25 \\
\kappa_P & 2 & \sigma_P & 0.4 \\
T & 1 & \rho & 0.6 \\
\bar{P} & 45 & \bar{X} & 12 \\
P_0 & 45 & X_0 & 15 \\ \hline
\end{array}$$
\end{minipage}&
\begin{minipage}{3in}
$$ \begin{array}{|cc|cc|} \hline
\multicolumn{2}{|c|}{\text{Producer 1}} & \multicolumn{2}{|c|}{\text{Producer 2}} \\ \hline
a_1 & 1  & a_2 & 2\\
b_1 &  2  & b_2 & 1\\
c_1 &  10  & c_2 & 80\\
g_1 & 8 & g_2 & 4 \\
K_1 &  0.2 & K_2 & 0.2 \\ \hline
\end{array}$$
\vspace{12pt}
\end{minipage}
\end{tabular*}
\caption{Model Parameters for the Examples in Section \ref{sec:numerical-ex}. \label{table:params}}
\end{table}

%
%
A large variety of CEPs are possible in our model; Table \ref{table:comp-corr} shows the game values corresponding to four representative correlation laws. These values were obtained by running Algorithm \ref{alg:main} discussed in Section \ref{sec:rmc} using $N=40000$ paths, and the basis functions $\{1, p, x, x^2, (2p - x- 80)_+, (p-2x - 10)_+ \}$. We find that different correlation laws modify the expected profit of the producers by $3\%-5\%$. As expected, individual producer values are maximized by the preferential equilibria that always favor the respective player. Counter-intuitively, the egalitarian CEP (which maximizes at each stage the minimum continuation value) produces larger game values to both producers than the utilitarian CEP (which maximizes the sum of continuation values). This occurs because the correlation law is applied \emph{stage-wise} and optimizes a local criterion; there is no guarantee that the corresponding global criterion is respected. A similar phenomenon was observed in \cite[Section 5.4]{RamseySzajowski08}.

\begin{table}
\begin{tabular}{|c|cc|}
\hline Correlation Law & $V_1(0,P_0,X_0)$ & $V_2(0, P_0, X_0)$ \\
\hline  Utilitarian & 5.30 & 4.14 \\
Egalitarian & 5.33 & 4.20 \\
Preferential 1 & 5.39 & 4.11\\
Preferential 2 & 5.02 & 4.24\\
 \hline
\end{tabular}
\medskip
\caption{Comparison of equilibrium game values for different correlation laws $\Gamma$. Standard errors of the Monte Carlo scheme are about 1\%. Parameters are as given in Table \ref{table:params}.\label{table:comp-corr}}
\end{table}

\begin{figure}[h]
\vspace*{-0.1in}\hspace*{-0.7in}{\includegraphics[width=7.8in]{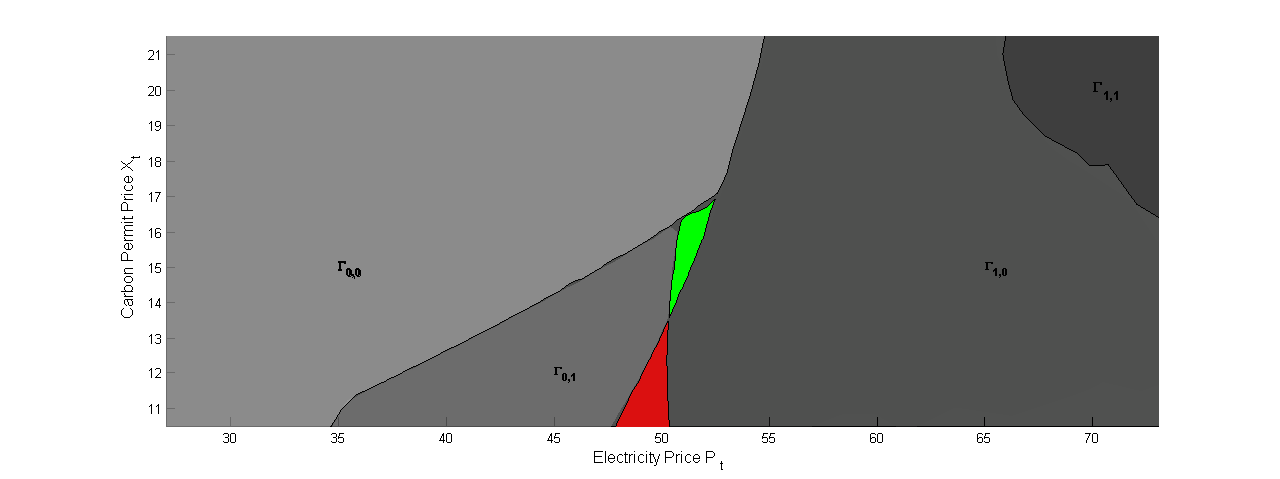}}
\caption{Equilibrium game strategy $\vxi^*(t)$ as a function of $(P_t, X^*_t)$ for $t=7$ and $\vze = (0,0)$. The green region denotes the anti-coordination game-type where the Preferential-1 correlation law is used, and the red region denotes the competitive game-type where the unique mixed NEP is chosen. Elsewhere, we label the regions according to the unique pure NEP implemented. \label{fig:two-producers-reg}}
\end{figure}

To illustrate the equilibrium strategy profiles, Figure \ref{fig:two-producers-reg} shows the empirical regions in the $(P,X)$-space corresponding to different equilibrium strategies at a fixed date $t=7$ (i.e.\ about three months into the year) using the Preferential-1 correlation law that always favors producer 1. As expected, when the current P\&L of both producers is strongly negative (upper-left corner), the equilibrium action is $\vxi^*(t) = (0,0)$; when it is strongly positive (large $P_t$) the equilibrium is to generate electricity $\vxi^*(t) = (1,1)$. Because of the differing carbon-efficiencies of the producers, there are also large regions where exactly one producer can generate profit (e.g.\ around $\{P_t \in [40, 45], X^*_t \in [10,12]\}$ only producer 2 is profitable). However at the border regions, the price impact and competition create new effects. In Figure \ref{fig:two-producers-reg}, we observe the emergence of a local anti-coordination game around $\{(P_t, X^*_t) = (50,15)\}$, and a competitive game around $\{(P_t, X^*_t) = (50,12)\}$.  We cannot analytically verify whether a particular type of game may emerge locally; thus the competitive game region in Figure \ref{fig:two-producers-reg} could be either a true phenomenon or an aberration due to numerical errors (e.g.\ poor regression fit in that region). Note that most simulated equilibrium paths for $X^*_t$ stay above $x=13$, so the competitive game scenario at $t=7$ is very unlikely to be realized (i.e.\ very few paths hit that region).

\begin{figure}[h]
\center{\includegraphics[width=6.7in]{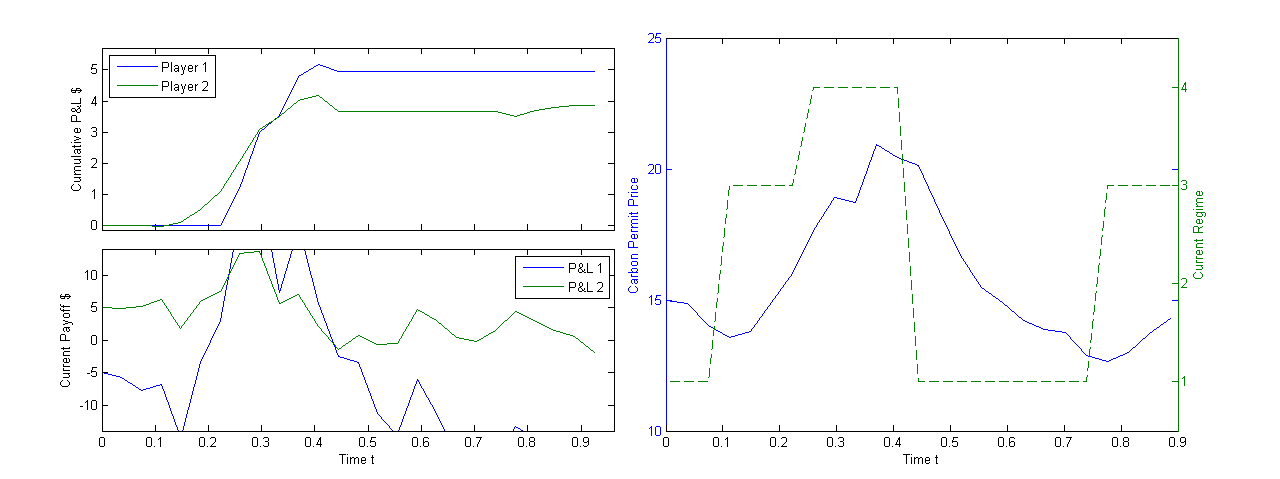}}
\caption{Sample equilibrium path of the emissions game.
Top left panel: cumulative realized P\&L of the players as a function of $t$. Bottom left panel: the electricity-carbon spread of each producer for the current time step. Right panel: evolution of the controlled equilibrium allowance price $X^*_t$, as well as the implemented strategy  $\vxi^*(t) \in \{00, 01, 10, 11\} \equiv \{1,\ldots, 4\}$. The panels were generated using Algorithm \ref{alg:forward-path} given in the Appendix. \label{fig:two-producers-path}}
\end{figure}

To better illustrate the optimal strategy over time, Figure \ref{fig:two-producers-path} shows a sample path of the equilibrium price $(X^*_t)$ for one $\omega$. Analogously to single-player problems, the $CO_2$ allowance price undergoes hysteresis cycles \cite{DixitPindyck}. Thus, when $(X^*_t)$ is low, production becomes profitable. This leads to increased emissions and $X^*_t$ tends to rise through the price impact mechanism. In turn, the ensuing higher emission costs eventually curtail production and $X^*_t$ falls back. The presence of switching costs $K_i$ lowers the scheduling flexibility of the producers and further amplifies this cycle through inertia.

\section{Conclusion}\label{sec:extensions}
In this paper we studied a new type of stochastic games which were motivated by dynamic emission schedules of energy producers under cap-and-trade schemes. Because multiple game equilibria can emerge, we explored various correlated equilibria. It is an interesting economic policy question which equilibrium is likely/desirable to be implemented and how the regulator can steer market participants towards that choice. For example, putting a price on emissions is supposed to partially drive out ``dirty'' producers. It would be an intriguing exercise to study how much these effects depend on equilibrium selection and whether \emph{blockading} of inefficient polluters is possible under some equilibria.

In our simplified model, the producers only made binary emission decisions at each stage. On a practical level, much finer granularity is available. It would be straightforward to extend our problem and allow a more general finite-state control set of size $| \mathcal{A}|$. The only modification would be to replace the $2 \times 2$ bimatrix games with a more general $\mathcal{A} \times \mathcal{A}$ bimatrix.
%
The theory for more than two producers is incomplete and it is an open problem to establish existence of CEP/NEP for multi-player stopping games (see \cite{SolanVieille02} for current state-of-the-art). 

\subsection{Further Extensions}\label{sec:cont-time}
Several aspects of our model merit further analysis. The dynamics for $C O_2$ allowance prices in \eqref{eq:X-dynamics} were selected to capture succinctly the price impact of each producer, leaving out other important features. As described in the introduction, as the permit expiration date $T$ approaches, the $CO_2$ price should converge either to zero (if excess permits remain) or to a fixed upper bound $\overline{x}$ (the penalty for emitting without an allowance). New (time-dependent) stochastic models are needed to mimic this property, see \cite{Carmona-Fehr,FehrHinz08}. Also, some cap-and-trade proposals will allow trading of allowances by financial participants whence \emph{no-arbitrage} restrictions might have to be imposed on the dynamics of $X$. All these possibilities can be handled straightforwardly, since the main construction is for arbitrary $X$-dynamics. Ideally, a fully endogenous model is desired for allowance prices; namely $X_t$ should be a function of total expected emissions until $T$ compared to total current supply, i.e.\ have a characterization in terms of conditional expectations of future equilibrium emission schedules. See \cite{CarmonaFehrHinz,Carmona-Fehr,ChesneyTaschini09} for such price-formation models and related
general equilibrium frameworks. These extensions will be considered in forthcoming papers.

Our formulation was in discrete-time; while this is sufficient for practical purposes, it is of great theoretical interest to construct a continuous-time model counterpart. The overall structure of a switching game as a sequence of stopping games straightforwardly carries over to continuous-time. However, description of correlated stopping equilibria in continuous time has not been attempted so far. In fact, the only reference dealing with randomized continuous-time stopping games is \cite{TouziVieille02} (see also \cite{LarakiSolanVieille05} for the latest results on general continuous timing games). Note that in continuous-time one must work with Nash $\eps$-equilibria since all stopping strategies are defined only in the almost-sure sense.  Second, to ensure the representation of $V_i$ as iterative stopping games through $\hJ^{n,m}_i$, it is necessary to \emph{a priori} show that each player makes finitely many regime switches. At this point we are not able to state any conditions to guarantee this, except requiring mandatory ``cool-off'' periods between each emission regime switch.

In our Markovian setting, solutions of continuous-time single-player switching problems have representations in terms of reflected backward stochastic differential equations (BSDE) \cite{HamadeneJeanblanc}. This representation should continue to hold in a game setting and will be explored in a separate paper. Related results have already been obtained for stochastic differential game analogues of our setup, whence $\vxi(t)$ has continuous state-space, see \cite{HamadeneLepeltierPeng97,HamadeneLepeltier00}.


\subsection*{Acknowledgment}
I am grateful to the anonymous referees for their valuable suggestions that improved the final presentation. Thanks also to Ren{\'e} Carmona for many useful discussions and David Ramsey, Krzysztof Szajowski, Jianfeng Zhang and participants at the MSRI Workshop on Economic Games and Mechanisms to Address Climate Change (May 2009) and the IPAM New Directions in Financial Mathematics Workshop (January 2010) for their feedback on earlier versions of this paper.

\bibliography{emissionSwitching}
\bibliographystyle{abbrv}

\appendix
\section*{Appendix: Numerical Algorithms}

\begin{algorithm}[b]
\caption{Simulating one realized cashflow path $\vth_i({s})$, $0 \le s \le T$ \label{alg:forward-path}}
\begin{algorithmic}
\REQUIRE Basis functions $B_\ell(p,x)$, $\ell = 1,\ldots,r$, regression coefficients $\valpha_i(t,\vze)$; correlation law $\Gamma$
\REQUIRE Initial condition $(p_0, x_0, \vxi(0))$; horizon $T$
\STATE   Initialize $\vth_i(0) \leftarrow 0 \quad$ \COMMENT{ Realized cashflows}
\FOR{$t=0,\ldots, T-1$}
     \FOR{ each $\vze \in \{0,1\}^2$}
        \STATE\COMMENT{Evaluate the predicted continuation values from taking action $\vze$}
        \STATE Set $ \hq_i(t, \vze) \leftarrow \sum_{\ell=1}^{r} \alpha^{\ell}_i(t, \vze) B_\ell( p_t, x_t) -K_{\{i,u_i(t),\zeta_i\}} + (a_i p_t - b_i x_t - c_i) \zeta_i$
     \ENDFOR
        \STATE Compute the stage-$t$ game values based on $\hq_{i}(t,\cdot)$, $i=1,2$  and $\Gamma$, see \eqref{eq:val-V-recursion}
        \STATE Obtain the correlated equilibrium strategy $\mathbf{\vxi}(t)$.
        \IF{ $\mathbf{\vxi}(t)$ is mixed}
          \STATE Perform randomization to obtain the realized action pair $\vxi({t+1})$
          \ELSE
          \STATE Set $\vxi({t+1}) \leftarrow \mathbf{\vxi}(t) \quad$ \COMMENT{ $\mathbf{\vxi}(t)$ is pure}
        \ENDIF
        \STATE Update $\vth_i({t+1}) \leftarrow \vth_i(t) -K_{\{i,u_i(t),u_i({t+1})\}} + (a_i p_t - b_i x_t - c_i) u_i({t+1}), \quad i=1,2$
        \STATE Make an independent draw  $(p_{t+1}, x_{t+1})\sim \PP^{\vxi({t+1})}(\cdot | p_t, x_t)$
\ENDFOR
\end{algorithmic}
\end{algorithm}

\begin{algorithm}[t]
\caption{Computing Correlated Equilibrium Game Values \label{alg:main}}
\begin{algorithmic}

\REQUIRE $N>0$ (number of paths); $B_\ell(p,x)$, $\ell=1,\ldots, r$ ($r$ regression basis functions)

\REQUIRE Correlation law $\Gamma$

\STATE Select anterior strategy profile $\vxi^0$

\FOR{each regime $\vze \in \{0,1\}^2$}

\STATE Simulate $N$ i.i.d.~paths $(p^n_t, x^{\vze,n}_t)_{n=1}^N$ under $\PP^{\vxi^0}$ using Algorithm \ref{alg:forward-path} and $p^n_0 = p_0$, $x^{\vze,n}_0 = x_0$
\ENDFOR

\STATE Initialize $\vth^{n}_i(T, \vze) \leftarrow 0$, $n=1,\ldots,N$

\FOR{ $t=(T-1), \ldots, 1, 0$}

 \FOR{ each regime $\vze$}

   \STATE Evaluate  $B_\ell(p^n_t, x^{\vze,n}_t)$ for $\ell=1,\ldots, r$ and $n=1,\ldots,N$

   \STATE Regress $$\valpha_i(t, \vze) \leftarrow \argmin_{\valpha \in \R^{r}}\sum_{n=1}^N  \Bigr|\vth^{n}_{i}(t+1, \vze) - \sum_{\ell=1}^{r} \alpha^\ell B_\ell(  p^n_t, x^{\vze,n}_t) \Bigl|^2$$

 \ENDFOR

  \FOR{ each current regime $\vxi$}
     \FOR{ each $\vze \in \{0,1\}^2$, and each $n=1,\ldots,N$}
        \STATE\COMMENT{Compute the predicted continuation value for each player from taking action $\vze$}
        \STATE Set $ \hq^{n}_i(t, \vxi, \vze) \leftarrow \sum_{\ell=1}^{r} \alpha^{\ell}_i(t, \vze) B_\ell( p^n_t, x^{\vxi,n}_t) -K_{\{i,u_i,\zeta_i\}} + (a_i p^n_t - b_i x^{\vxi,n}_t - c_i) \zeta_i$.
     \ENDFOR
     \FOR{ each path $n=1,\ldots, N$}
        \STATE Compute the stage-$t$ game values based on $\hq^{n}_{\cdot}(t,\vxi,\cdot)$  and $\Gamma$, see \eqref{eq:val-V-recursion}
        \STATE Obtain the equilibrium policy ${\vxi}^{n,*}(t,\vxi)$
        \STATE Recompute $\vth^n_i(t, \vxi)$ using ${\vxi}^{n,*}(t,\vxi)$ at stage-$t$ and Algorithm \ref{alg:forward-path} for future stages
     \ENDFOR
  \ENDFOR

\ENDFOR

  \RETURN $V_i(0, p_0, x_0, \vze) \simeq \frac{1}{N} \sum_{n=1}^N \vth^{n}_i(0, \vze)$
  \RETURN Regression coefficients $\valpha_i(t,\vze)$ summarizing equilibrium strategies
  \end{algorithmic}
 \end{algorithm}

\end{document}